\documentclass[12pt]{article}
\usepackage{amsmath,amsthm,amscd,amsfonts}

\newcounter{num}
\newtheorem{theorem}{Theorem}[section]

\theoremstyle{plain}

\newtheorem*{Moore}{Moore's Conjecture}

\newtheorem*{Serre}{Serre's Theorem}

\newtheorem*{main}{Theorem}

\theoremstyle{remark}

\theoremstyle{definition}
\newtheorem{definition}[theorem]{Definition}

 \DeclareMathOperator{\Ext}{Ext}

\begin{document}

\newcommand{\be}        {\begin{eqnarray}}
\newcommand{\ee}        {\end{eqnarray}}
\newcommand{\pl}{\partial}
\newcommand{\sbs}{\subset}
\newcommand{\vr}{\varphi}
\newcommand{\lm}{\lambda}
\newcommand{\eps}{\varepsilon}
\newcommand{\nb}{\nabla}
\newcommand{\wt}{\widetilde}

\newcommand{\cV}      {{\cal V}}
\newcommand{\cE}      {{\cal E}}
\newcommand{\cD}      {{\cal D}}
\newcommand{\cK}      {{\cal K}}
\newcommand{\cB}      {{\cal B}}
\newcommand{\cI}      {{\cal I}}
\newcommand{\cM}      {{\cal M}}
\newcommand{\cA}      {{\cal A}}
\newcommand{\cR}      {{\cal R}}
\newcommand{\cT}      {{\cal T}}
\newcommand{\cP}      {{\cal P}}
\newcommand{\cC}      {{\cal C}}
\newcommand{\cQ}      {{\cal Q}}
\newcommand{\cG}      {{\cal G}}
\newcommand{\cW}      {{\cal W}}
\newcommand{\cL}      {{\cal L}}
\newcommand{\cF}      {{\cal F}}
\newcommand{\cS}      {{\cal S}}
\newcommand{\cH}      {{\cal H}}
\newcommand{\cN}      {{\cal N}}
\newcommand{\AAA}       {{\Bbb A}}     
\newcommand{\aaa}       {{\LBbb A}}
\newcommand{\BB}        {{\Bbb B}}
\newcommand{\bb}        {{\LBbb B}}
\newcommand{\CC}        {\mathbb{C}}
\newcommand{\cc}        {{\LBbb C}}
\newcommand{\DD}        {{\Bbb D}}
\newcommand{\dd}        {{\LBbb D}}
\newcommand{\EEE}       {{\Bbb E}}
\newcommand{\eee}       {{\LBbb E}}
\newcommand{\FF}        {{\Bbb F}}
\newcommand{\ff}        {{\LBbb F}}
\newcommand{\GGG}       {{\Bbb G}}
\newcommand{\ggg}       {{\LBbb G}}
\newcommand{\HH}        {{\Bbb H}}
\newcommand{\hh}        {{\LBbb H}}
\newcommand{\II}        {{\Bbb I}}
\newcommand{\ii}        {{\LBbb I}}
\newcommand{\JJ}        {{\Bbb J}}
\newcommand{\jj}        {{\LBbb J}}
\newcommand{\KK}        {{\Bbb K}}
\newcommand{\kk}        {{\LBbb k}}
\newcommand{\LLL}       {{\Bbb L}}
\newcommand{\lll}       {{\LBbb L}}
\newcommand{\MM}        {{\Bbb M}}
\newcommand{\mm}        {{\LBbb M}}
\newcommand{\NN}        {\mathbb{N}}
\newcommand{\nn}        {{\LBbb N}}
\newcommand{\OO}        {{\Bbb O}}
\newcommand{\oo}        {{\LBbb O}}
\newcommand{\PP}        {{\Bbb P}}
\newcommand{\pp}        {{\LBbb P}}
\newcommand{\QQ}        {{\Bbb Q}}
\newcommand{\qq}        {{\LBbb Q}}
\newcommand{\RR}       { \mathbb{R}}
\newcommand{\rr}        {{\LBbb R}}
\newcommand{\SSS}       {{\Bbb S}}
\newcommand{\sss}       {{\LBbb S}}
\newcommand{\TTT}       {{\Bbb T}}
\newcommand{\ttt}       {{\LBbb t}}
\newcommand{\UU}        {{\Bbb U}}
\newcommand{\VV}        {{\Bbb V}}
\newcommand{\vv}        {{\LBbb V}}
\newcommand{\WW}        {{\Bbb W}}
\newcommand{\ww}        {{\LBbb W}}
\newcommand{\XX}        {{\Bbb X}}
\newcommand{\xx}        {{\LBbb X}}
\newcommand{\YY}        {{\Bbb Y}}
\newcommand{\yy}        {{\LBbb Y}}
\newcommand{\ZZ}        {\mathbb{Z}}
\newcommand{\zz}        {{\LBbb Z}}

\title{On Cohomology rings of infinite groups}
\author{{Eli Aljadeff}\thanks{Supported by E. and J. Bishop Research Fund.}\\
Department of Mathematics\\
Technion -- Israel Institute of Technology\\
32000 Haifa, Israel\\
{\small e-mail: aljadeff@math.technion.ac.il;
aljadeff@tx.technion.ac.il}}



\maketitle

\vspace{1cm}

\begin{abstract}

Let $R$ be any ring (with $1)$, $G $ a torsion free group and $RG
$ the corresponding group ring. Let $Ext_{RG}^{*}(M,M)$ be the
cohomology ring associated to the $RG $-module $M$. Let $H$ be a
subgroup of finite index of $G $. The following is a special
version of our main Theorem: Assume the profinite completion of $G
$ is torsion free. Then an element $\zeta $ $\in$
$\Ext_{RG}^{*}(M,M)$ is nilpotent (under Yoneda's product) if and
only if its restriction to $\Ext_{RH}^{*}(M,M)$ is nilpotent. In
particular this holds for the Thompson group.

There are torsion free groups for which the analogous statement is
false.

\end{abstract}

\vspace{1cm}

\section{Introduction.}

Elementary abelian subgroup induction plays a crucial role in
cohomology and representation theory of finite groups (see
\cite{Se}, \cite{Qu1}, \cite{Qu2}, \cite{QV}, \cite{Ca1},
\cite{Ca2}, \cite{Ca3}, \cite{Ch}, \cite{AS}, \cite{Ad},
\cite{AM}, \cite{Be1}, \cite{AE}, \cite{AG}). Roughly speaking,
the results say that important cohomological properties hold for a
group ring $RG$, $G$ finite and $R$ an arbitrary ring (with {1}),
if and only if they hold for $RE$ where $E$ runs over all
elementary abelian subgroups of $G$. In general, similar
statements are false if one replaces the family of elementary
abelian subgroups by cyclics.

We recall two important results which exhibit the role of the
elementary abelian subgroups:

\begin{main}[Quillen and Carlson {\cite{Qu1}, \cite{Ca2}, \cite{Ca3}}]
Let $G $ be a finite group and $R$ an arbitrary ring. Let $M$ be a
module over $RG$ and let $\Ext_{RG}^{*}(M,M)$ be the associated
cohomology ring (with Yoneda's product). Then an element $\zeta$
$\in$ $\Ext_{RG}^{*}(M,M)$ is nilpotent if and only if its
restriction to $\Ext_{RE}^{*}(M,M)$ is nilpotent where $E$ runs
over all elementary abelian subgroups of $G$.

\end{main}

\begin{main}[Chouinard \cite{Ch}]
Let $G $ be a finite group and $R$ be an arbitrary ring. If $M$ is
any module over $RG $ then it is weakly projective (projective) if
and only if it is weakly projective (projective) over all subrings
$RE$ where $E$ is an elementary abelian subgroup of $G.$

\end{main}

Our main objective in this paper is to prove an ``infinite version'' of Quillen and Carlson
Theorem.

In 1976 Moore posed the following conjecture which may be viewed as an ``infinite version'' of
Chouinard's theorem.

\begin{Moore}[see \cite{Ch}, {\cite[Conjecture~1.1]{ACGK}}]\label{Moore}
  Let $G $ be a group and $H $ a subgroup of finite index. Let
$R$ be an arbitrary ring. Assume that for every nontrivial element
$x$ in $G $, at least one of the following two conditions holds:

M1)~$\langle x\rangle \cap H \neq \{e\}$ (in particular this holds
if $G $ is torsion free)

M2)~$ord(x)$ is finite and invertible in $R.$ \\
Then every $RG -$module $M$ which is projective over $RH,$ is
projective also over $RG .$
\end{Moore}
We refer to (M1) or (M2) as Moore's condition for the triple $(G,
H, R)$.

Note that Chouinard's theorem implies Moore's conjecture whenever
the group $G $ is finite.

Moore's conjecture is a far reaching generalization of the
following result of Serre which is well know.

\begin{Serre}[see \cite{Sw}]\label{Serre}
  Let $G $ be a group and $H$ a subgroup of finite
index. Assume that $H$ has finite cohomological dimension (that is
$G $ has virtual finite cohomological dimension). If $G $ is
torsion free then it has finite cohomological dimension. Moreover,
$cd(G )=cd(H).$
\end{Serre}

Let us show how Serre's Theorem is obtained from Moore's
conjecture. Assume $cd(H)=n$. If $P.$ $\rightarrow
{\ZZ}\rightarrow 0$ is a projective resolution of ${\ZZ}$ over
${\ZZ}G ,$ it is projective also over $H.$ It follows that the
$n$-th syzygy $Y_{n}$ of the resolution is a ${\ZZ}G $- module
whose restriction to $H$ is projective. Moore's conjecture says
that $ Y_{n}$ is projective over ${\ZZ}G $ and so $cd(G )\leq n.$

Our main goal is to formulate and prove a ``Moore's analog'' to
Quillen and Carlson's result. We say "an analog" and not ``the
analog'' since ``the analog'' is false in its full generality as
observed by Dave Benson (see Remark 2 below). Before stating the
theorem precisely we set some terminology and notation. Let $G $
be any group and denote by $P(G )$ the collection of all finite
index, normal subgroups of $ G.$ Let $\Omega$=$\Omega _{G }$ be a
subset of $P(G )$ filtered from below. Assume further that
$\Omega$ is cofinal in $P(G )$ and let $\widehat{G
}_{\Omega}={\lim }G /N$ be the profinite completion of $G $ with
respect to $\Omega$. Let $\phi :G \rightarrow \widehat{G
}_{\Omega}$ be the canonical map induced by the natural
projections $G \rightarrow G /N,\;N\in \Omega$ (see \cite{RZ}).


\begin{definition}\label{D1.1}
We say that $\widehat{G }_{\Omega}$ has no new torsion (with
respect to $G $) if any element $z$ of prime order (say $p$) in
$\widehat{G }_{\Omega}$ is conjugate to an element $\phi(x)$ where
$x$ is an element (in $G$) of order $p$.
\end{definition}

 \noindent
{\it Remark (1).}~~(i)~if $\widehat{G }_{\Omega}$ is torsion free
then clearly it has no new torsion. (ii) no assumption has been
made about $G $ being residually finite.

\begin{theorem}\label{Th1.2}
Let $RG$ be a group ring over $R$. Let $H$ be a subgroup of finite
index in $G$. Assume Moore's condition (M1) or (M2) holds. If
$\widehat{G }_{\Omega}$ has no new torsion ($\Omega$ as above)
then  an element $\zeta$ in $\Ext_{RG}^{*}(M,M)$ is nilpotent if
and only if its restriction to $\Ext_{RH}^{*}(M,M)$ is nilpotent.
\end{theorem}

In particular this holds for the Thompson group $T$ (since
$\widehat{T}_{P(T)}$ is torsion free). (See \cite{BG} for
the definition of the Thompson group and some of its properties).

  Quillen and Carlson's result was generalized to infinite
groups $G $ which have virtual finite cohomological dimension. The
idea already appears in Quillen's paper (\cite{Qu2}). The result
says that a nilpotent element in the cohomology ring is detected
by its restrictions to the elementary abelian subgroups. This
result was extended considerable by Benson, namely to groups which
belong to Kropholler's class ${\scriptstyle\bf LH}\frak F$
(Locally-${\scriptstyle\bf H}\frak F$) and modules which are
$FP_\infty$ (see \cite{Be2}, \cite{Be3}). It is important to
mention that the Thompson group does not belong to the class
${\scriptstyle\bf LH}\frak F$ (see \cite{Kr}).

\vspace{0.5 cm}

 \noindent
{\it Remark (2).}~~Dave Benson observed that if one drops the
assumption on the profinite completion of $G$ then the analogous
statement to Theorem~\ref{Th1.2} is false in general. Indeed,
based on the Kan-Thurston construction in \cite{KT}, it is shown
in \cite{BDH}, that there exists a torsion free group $G$ whose
mod-$p$ cohomology is isomorphic to the mod-$p$ cohomology of
$Z/p$ (the cyclic group of order $p$). Furthermore, $G$ has a
perfect subgroup $P$ of index $p$ whose mod-$p$ cohomology is
trivial in positive degrees (see also [CL], where it is shown that
$G$ does not satisfy the Quillen-Carlson theorem).

\section{Proofs.}

An important tool in the proof of Theorem~\ref{Th1.2} is the
crossed product construction. It allows us to use induction,
namely to express the crossed product $R*G$ in terms of a
subalgebra $R*H$ and $G/H$ where $H$ is a normal subgroup of $G$.
Obviously this is not possible within the family of groups rings.
Therefore we shall work in the context of crossed products rather
than group rings. The first step is to extend Quillen--Carlson's
result to arbitrary crossed products $R*G$ where $G$ is finite.

\begin{theorem}\label{Th2.1}
Let $M$ be a module over a crossed product $R*G$, where $R$ is any
ring and $G$ a finite group. Then an element $\zeta$ $\in$
$\Ext_{R*G}^{*}(M,M)$ is nilpotent if and only if its restriction
to $\Ext_{R*E}^{*}(M,M)$ is nilpotent where $E$ runs over all
elementary abelian subgroups of $G$.
\end{theorem}

\begin{proof}
For the proof we need to generalize Corollary~2.2, Theorem~2.5 and
Lemma~2.6 in \cite{Ca3} to arbitrary crossed products $R*G$. The
generalizations are straightforward and left to reader.
\end{proof}

Next, we extend Theorem~\ref{Th1.2}. This is a two step extension. For clarity we do one at a time.
The first one is to the context of crossed products.

\begin{theorem}\label{Th2.2}
 Let $R*G$ be a crossed product of the group $G$ over an arbitrary ring $R$. Let $H$ be
a subgroup of finite index in $G$. Let $M$ be an $R*G$ module.
Assume Moore's condition (M1) or (M2) holds. If $\widehat{G
}_{\Omega}$ has no new torsion then an element $\zeta$ $\in$
$\Ext_{R*G}^{*}(M,M)$ is nilpotent if and only if its restriction
to $\Ext_{R*H}^{*}(M,M)$ is nilpotent.
\end{theorem}

The second extension consists in relaxing the condition on the
torsion in $\widehat{G }_{\Omega}$. It is reasonable to expect
that one should be careful only about the torsion outside the
completion of $H$. In order to state this condition precisely let
$G $ be any group and $H$ a subgroup of finite index. Let
$\Omega=\Omega_{G, H} $ be a collection of subgroups of $H,$
filtered from below, normal and of finite index in $G .$ Let
$\widehat{G }_{\Omega}$ and $\widehat{H }_{\Omega}$ be the
corresponding profinite completions of $G $ and $H$.


\begin{definition}\label{D2.3}
 We say that $\widehat{G }_{\Omega}$ has no new torsion outside $\widehat{H }_{\Omega}$ if
any element of prime order (say $p)$ in $\widehat{G
}_{\Omega}\,\backslash \,\widehat{H }_{\Omega}$ is conjugate to an
element $\phi (x)$ where $ x$ is an element in $G $ of order $p.$
\end{definition}

\begin{theorem}\label{Th2.4}
 Let $R*G$ be a crossed product of the group $G$ over an arbitrary ring $R$. Let $H$ be
a subgroup of finite index in $G$. Let $M$ be an $R*G$ module.
Assume Moore's condition (M1) or (M2) holds. If $\widehat{G
}_{\Omega}$ ($\Omega$=$\Omega_{G, H} $) has no new torsion outside
$\widehat{H }_{\Omega}$ then an element $\zeta$ $\in$
$\Ext_{R*G}^{*}(M,M)$ is nilpotent if and only if its restriction
to $\Ext_{R*H}^{*}(M,M)$ is nilpotent.
\end{theorem}

The similarity of Theorem~\ref{Th2.4} above and Theorem 2.2 in [Al] leads us to formulate a more
general result.

Let $P$ be a property associated to an arbitrary crossed product
$R*G$ which satisfies the following conditions:
\begin{enumerate}
\item[1)] If $P$ holds for $R*G$ then it holds for $R*H$ where $H$
is any subgroup of $G$.

\item[2)] If $G$ is finite, then $P$ holds for $R*G$ if and only
if it holds for $R*E$ where $E$ runs over all elementary abelian
subgroups of $G$.
\end{enumerate}

Of course, the properties $P$ considered in this article and in
\cite{Al} are: 1) for an $R*G$ module $M$, an element $\zeta$
$\in$ $\Ext_{RG}^{*}(M,M)$ is nilpotent (under Yoneda's product);
2) an $R*G$ module $M$ is weakly projective (projective).

For such $P$ we have

\begin{theorem}\label{Th2.5}
Let $R$ be any ring and let $G $ be any group. Let $H$ be a
subgroup of finite index of $G$ which contains all elements of
order $p$ (prime) for $p$ not invertible in $R$ (i.e. (M1) or (M2)
holds). If $\widehat{G }_{\Omega}$ ($\Omega$=$\Omega_{G, H} $) has
no new torsion outside $\widehat{H }_{\Omega}$, then property $P$
holds for $R*G$ if and only if it holds for $R*H$.

\end{theorem}

\begin{proof}
The proof is identical to the proof of Theorem~3.1 in \cite {Al} and also is left to the reader.
\end{proof}

Acknowledgment: I would like to thank Dave Benson for his important comment with regard the
main result of the article (see Remark 2).

\end{document}